\documentclass[11pt]{amsart}
\usepackage{amsmath}
\usepackage{amssymb}
\usepackage{latexsym}
\usepackage{graphicx}

\newcommand{\N}{\mathbb{N}}

\newtheorem{theorem}{Theorem}
\newtheorem{corollary}{Corollary}
\newtheorem{lemma}{Lemma}
\theoremstyle{remark}
\newtheorem{remark}{Remark}

\begin{document}

\title[$n$-starlike integral operators]{\large On some $n$-starlike integral operators}

\author[K. O. Babalola]{K. O. BABALOLA}

\begin{abstract}
By a completion of a lemma of Babalola and Opoola \cite{KT}, we prove that certain generalized integral operators preserve $n$-starlikeness in the open unit disk $E=\{z\in \mathbb{C}: |z|<1\}$. Our results generalize, extend and improve many known ones.
\end{abstract}



\maketitle

\section{Introduction}
Let $A$ be the class of functions
\begin{equation}
f(z)=z+a_2z^2+\cdots\, \label{1}
\end{equation}
which are analytic in $E$. A function $f\in A$ is said to be starlike of order $\lambda$, $0\leq\lambda<1$ if and only if, for $z\in E$,
\[
Re\;\frac{zf'(z)}{f(z)}>\lambda.
\]

Also a function $f\in A$ is said to be convex of order $\lambda$, $0\leq\lambda<1$ if and only if, for $z\in E$,
\[
Re\;\left\{1+\frac{zf''(z)}{f'(z)}\right\}>\lambda.
\]

Let $S^\ast(\lambda)$ and $K(\lambda)$ denote, as usual, the classes of starlike and convex functions of order $\lambda$ respectively. Salagean \cite{GS} introduced the operator $D^n$, $n\in\N$ as:
\[
D^nf(z)=D(D^{n-1}f(z))=z[D^{n-1}f(z)]' 	 
\]
with $D^0f(z)=f(z)$ and used it to generalize the concepts of starlikeness and convexity of functions in the unit disk as follows: a function $f\in A$ is said to belong to the classes $S_n(\lambda)$ if and only if
\[
Re\frac{D^{n+1}f(z)}{D^nf(z)}>\lambda.
\]
For $n=0,1$, we have the classes of starlike and convex functions respectively. We refer to functions of the classes $S_n(\lambda)$ as $n$-starlike functions in the unit disk. For $\lambda=0$ we simply write $S^\ast$, $K$ and $S_n$.

Let $\beta>0$, $\alpha\geq 0$ be real numbers, $\gamma$ and $\delta$ complex constants with $\alpha+\delta=\beta+\gamma$. For $f\in A$, the generalized integral operator
\begin{equation}
\mathcal{J}(f)=\left\{\frac{\beta+\gamma}{z^\gamma}\int_0^zt^{\delta-1}f(t)^\alpha dt\right\}^{\frac{1}{\beta}},\;\;\;\beta+Re\;\gamma\geq 0,\ \label{2}
\end{equation}
and its many special cases (for example: $\beta=\alpha=1$, $\gamma=\delta$; $\beta=\alpha=1$, $\gamma=\delta=0$; $\beta=\alpha=1$, $\gamma=1$ and $\delta=1-\alpha$) have been studied repeatedly in many literatures \cite{SA,BO,SD,WM,JK,ZL,RJ,EP,MM,MN,TO,WC,RS} where $f(z)$ belongs to some favoured classes of functions. More general integral operators were studied in \cite{MM} where the authors used a new method of analysis to obtain results that are both more general and sharper than many earlier ones. 

Let $\beta>0$, $\alpha\geq 0$ be real numbers, $\gamma$ and $\delta$ complex constants such that $\alpha+\delta=\beta+\gamma$. Define $\mathcal{J}_0^j(z)^\beta=f(z)^\alpha$, $j=1,2$ and for $m\in N$ define 
\[
\mathcal{J}_m^1(f)=\left\{\frac{(\beta+\gamma)^m}{z^\gamma\Gamma(m)}\int_0^z\left(\log\frac{z}{t}\right)^{m-1}t^{\delta-1}f(t)^\alpha dt\right\}^{\frac{1}{\beta}},
\]
where Re $\gamma\geq 0$ and
\[
\mathcal{J}_m^2(f)=\left\{\binom{\beta+\gamma+m-1}{\beta+\gamma-1}\frac{m}{z^\gamma}\int_0^z\left(1-\frac{t}{z}\right)^{m-1}t^{\delta-1}f(t)^\alpha dt\right\}^{\frac{1}{\beta}}
\]
also with $m-1$+Re $\gamma\geq 0$.

The integrals $\mathcal{J}^j(f)$ are similar to the Jung-Kim-Srivastava one-parameter families of integral operators \cite{JK}. However, only in the case $\beta=\alpha=1$ and $\gamma$ real, then $\mathcal{J}^j(f)$ are special cases of those in \cite{JK}. Furthermore if $m=1$, both integrals yield the integral operator ~(\ref{2}).

In the present paper, we will study the integrals $\mathcal{J}^j(f)$ for $f$ belonging to the classes $S_n(\lambda)$. Furthermore, if $\gamma$ and $\delta$ are real constants we will obtain the best possible inclussion for $\mathcal{J}^j(f)$ given that $f\in S_n(\lambda)$. Natural corollaries to the main results of this work are that: (i) for all real number $\beta$, $\alpha\geq 0$, the integrals $\mathcal{J}^j(f)$, $j=1,2$ preserve starlikeness and convexity in the open unit disk and that (ii) our result will improve and extend many known ones for all the many special cases. The main results are presented in Section 3 while we discuss the many special cases arising from taking $m=1$ in Section 4.

In the next section we give some lemmas necessary for the proof of our results.

\section{Preliminary Lemmas}
Let $P$ denote the class of functions $p(z)=1+c_1z+c_2z^2+\cdots$ which are regular in $E$ and satisfy Re $p(z)>0$, $z\in E$. We shall need the following lemmas.

\begin{lemma}[\cite{BO}]
Let $u=u_1+u_2i$, $v=v_1+v_2i$ and $\psi(u,v)$ a complex-valued function satisfying:
 
{\rm(a)} $\psi(u,v)$ is continuous in a domain $\Omega$ of $\mathbb{C}^2$,

{\rm(b)} $(1,0)\in\Omega$ and Re$\psi(1,0)>0$,

{\rm(c)} Re$\psi(\lambda+(1-\lambda)u_2i, v_1)\leq\lambda$ when $(\lambda+(1-\lambda)u_2i, v_1)\in\Omega$ and $2v_1\leq -(1-\lambda)(1+u_2^2)$ for real number $0\leq\lambda<1$.

If $p\in P$ such that $(p(z),zp'(z))\in\Omega$ and $Re$ $\psi(p(z),zp'(z))>\lambda$ for $z\in E$, then $Re$ $p(z)>\lambda$ in $E$.
\end{lemma}

The above lemma is an abridged form of a more detail one in \cite{BO}.

\begin{lemma}[\cite{EE}]
Let $\eta$ and $\mu$ be complex constants and $h(z)$ a convex univalent function in $E$ satisfying $h(0)=1$, and $Re(\eta h(z)+\mu)>0$. Suppose $p\in P$ satisfies the differential subordination:
\begin{equation}
p(z)+\frac{zp'(z)}{\eta p(z)+\mu}\prec h(z),\;\;\;z\in E.\, \label{3}
\end{equation}
If the differential equation:
\begin{equation}
q(z)+\frac{zq'(z)}{\eta q(z)+\mu}=h(z),\;\;\;q(0)=1\, \label{4}
\end{equation}
has univalent solution $q(z)$ in $E$, then $p(z)\prec q(z)\prec h(z)$ and $q(z)$ is the best dominant in $~(\ref{3})$.
\end{lemma}

The formal solution of ~(\ref{4}) is given as
\[
q(z)=\frac{zF'(z)}{F(z)}=\frac{\eta+\mu}{\eta}\left(\frac{H(z)}{F(z)}\right)^\eta-\frac{\mu}{\eta}
\]
where
\[
F(z)^\eta=\frac{\eta+\mu}{z^\mu}\int_0^zt^{\mu-1}H(t)^\eta dt
\]
and
\[
H(z)=z.\exp\left(\int_0^z\frac{h(t)-1}{t}dt\right)
\]
(see \cite{SS,HM}). The authors in \cite{SS} gave sufficient conditions for the univalence of the solution, $q(z)$, of ~(\ref{4}) as well as some generalised univalent solutions for some given $h(z)$. 

The second part of the next lemma is the completion of Lemma 2.2 in \cite{KT}.

\begin{lemma}[\cite{KT}]
Let $f\in A$ and $\zeta>0$ be real.

{\rm (i)} If for $z\in E$, $D^{n+1}f(z)^\zeta/D^nf(z)^\zeta$ is independent of $n$, then
\begin{equation}
\frac{D^{n+1}f(z)^\zeta}{D^nf(z)^\zeta}=\zeta\frac{D^{n+1}f(z)}{D^nf(z)}.\, \label{5}
\end{equation}

{\rm (ii)} The equality $~(\ref{5})$ also holds if $D^{n+1}f(z)/D^nf(z)$ is independent of $n$, $z\in E$.
\end{lemma}

\begin{proof}
The proof of the first part of the above lemma was presented in \cite{KT}. As for (ii), let $D^{n+1}f(z)/D^nf(z)$ assume the same value for all $n\in\N$. For $n=0$, the assertion is easy to verify. Let $n=1$. Then
\begin{align*}
\frac{D^2f(z)^\zeta}{D^1f(z)^\zeta}
&=1+\frac{zf''(z)}{f'(z)}+(\zeta-1)\frac{zf'(z)}{f(z)}\\
&=\frac{D^2f(z)}{D^1f(z)}+(\zeta-1)\frac{D^1f(z)}{D^0f(z)}.
\end{align*}
Since $D^1f(z)/D^0f(z)=D^2f(z)/D^1f(z)$ we have
\[
\frac{D^2f(z)^\zeta}{D^1f(z)^\zeta}=\zeta\frac{D^2f(z)}{D^1f(z)}.
\]
Now suppose ~(\ref{5}) holds for some integer $k$. Then
\begin{equation}
\frac{D^{k+2}f(z)^\zeta}{D^{k+1}f(z)^\zeta}=\frac{D^{k+2}f(z)}{D^{k+1}f(z)}+(\zeta-1)\frac{D^{k+1}f(z)}{D^kf(z)}.\, \label{6}
\end{equation}
Since $D^{k+1}f(z)/D^kf(z)$ has the same value for each $k\in\N$, we can write ~(\ref{6}) as
\[
\frac{D^{k+2}f(z)^\zeta}{D^{k+1}f(z)^\zeta}=\frac{D^{k+2}f(z)}{D^{k+1}f(z)}+(\zeta-1)\frac{D^{k+2}f(z)}{D^{k+1}f(z)}
\]
which implies
\[
\frac{D^{k+2}f(z)^\zeta}{D^{k+1}f(z)^\zeta}=\zeta\frac{D^{k+2}f(z)}{D^{k+1}f(z)}.
\]
Thus the lemma follows by induction.
\end{proof}

\begin{remark}
Let $f\in S_n(\lambda)$. Then there exists $p\in P$ such that
\[
\frac{D^{n+1}f(z)}{D^nf(z)}=\lambda+(1-\lambda)p(z)
\]
independent of $n\in \N$. Hence for $f\in S_n(\lambda)$, the assertion of Lemma 2 holds. Thus we have
\[
Re\frac{D^{n+1}f(z)^\zeta}{D^nf(z)^\zeta}=\zeta Re\frac{D^{n+1}f(z)}{D^nf(z)}>\zeta\lambda.
\]
In particular, if $\lambda=0$, then for $\zeta>0$ we have Re $\frac{D^{n+1}f(z)^\zeta}{D^nf(z)^\zeta}>0$ if and only if Re $\frac{D^{n+1}f(z)}{D^nf(z)}>0$.
\end{remark}

\section{Main Results}

\begin{theorem}
Let $\alpha\geq 0$. Suppose for $\alpha>0$, the real number $\lambda$ is defined such that $0\leq\alpha\lambda<1$. If $f\in S_n(\lambda)$, then $\mathcal{J}^j(f)\in S_n(\frac{\alpha}{\beta}\lambda)$, $j=1,2$.
\end{theorem}

\begin{proof}
Let $f\in S_n(\lambda)$ have the form ~(\ref{1}). If $\alpha=0$, then $\mathcal{J}^j(f)=z$ by evaluation using the Beta and Gamma functions. Thus the result holds trivially in this case. Suppose $\alpha>0$, then we can write
\[
f(z)^\alpha=z^\alpha+A_2(\alpha)z^{\alpha+1}+...
\]
where $A_k(\alpha)$, $k=2,3,...$, depends on the coefficients $a_k$ of $f(z)$ and the index $\alpha$. Thus evaluating the integrals in series form, also using the Beta and Gamma functions and noting that
\[
\binom{\sigma}{\gamma}=\frac{\Gamma(\sigma+1)}{\Gamma(\sigma-\gamma+1)\Gamma(\gamma+1)}
\]
we obtain
\[
\mathcal{J}_m^1(f)^\beta=z^\beta+\sum_{k=2}^\infty\left(\frac{\beta+\gamma}{\beta+\gamma+k-1}\right)^m A_k(\alpha)z^{\beta+k-1}
\]
and
\[
\mathcal{J}_m^2(f)^\beta=z^\beta+\frac{\Gamma(\beta+\gamma+m)}{\Gamma(\beta+\gamma)}\sum_{k=2}^\infty\frac{\Gamma(\beta+\gamma+k-1)}{\Gamma(\beta+\gamma+m+k-1)}A_k(\alpha)z^{\beta+k-1}.
\]
From the above series expansions we can see that $\mathcal{J}_0^j(f)^\beta=f(z)^\alpha$, $j=1,2$ are well defined. Also from the series expansions we find the recurssive relation
\begin{equation}
\mu \mathcal{J}_m^j(z)^\beta+z(\mathcal{J}_m^j(f)^\beta)'=\xi\mathcal{J}_{m-1}^j(f)^\beta,\;\;\;j=1,2\, \label{7}
\end{equation}
where $\mu=\gamma$ and $\xi=\beta+\gamma$ for $j=1$ and $\mu=\gamma+m-1$ and $\xi=\beta+\gamma+m-1$ if $j=2$. Furthermore let $\mu=\mu_1+\mu_2i$. Now applying the operator $D^n$ on ~(\ref{7}) we have
\[
\frac{D^{n+1}\mathcal{J}_{m-1}^j(f)^\beta}{D^n\mathcal{J}_{m-1}^j(f)^\beta}=\frac{\mu D^{n+1}\mathcal{J}_m^j(f)^\beta+D^{n+2}\mathcal{J}_m^j(f)^\beta}{\mu D^n\mathcal{J}_m^j(f)^\beta+D^{n+1}\mathcal{J}_m^j(f)^\beta}.
\]
Let $p(z)=\frac{D^{n+1}\mathcal{J}_m^j(z)^\beta}{D^n\mathcal{J}_m^j(z)^\beta}$. Then
\begin{equation}
\frac{D^{n+1}\mathcal{J}_{m-1}^j(z)^\beta}{D^n\mathcal{J}_{m-1}^j(z)^\beta}=p(z)+\frac{zp'(z)}{\mu+p(z)}.\, \label{8}
\end{equation}
Define $\psi(p(z),zp'(z))=p(z)+\frac{zp'(z)}{\mu+p(z)}$ for $\Omega=[\mathbb{C}-\{-\mu\}]\times\mathbb{C}$. Obviously $\psi$ satisfies the conditions (a) and (b) of Lemma 1. Now let $0\leq\lambda_0=\alpha\lambda<1$. Then $\psi(\lambda_0+(1-\lambda_0)u_2i,v_1)=\lambda_0+(1-\lambda_0)u_2i+\tfrac{v_1}{\mu+(\lambda_0+(1-\lambda_0)u_2i)}$ so that Re $\psi(\lambda_0+(1-\lambda_0)u_2i,v_1)=\lambda_0+\tfrac{(\mu_1+\lambda_0)v_1}{(\mu_1+\lambda_0)^2+(\mu_2+(1-\lambda_0)u_2)^2}$. If $v_1\leq-\tfrac{1}{2}(1-\lambda_0)(1+u_2^2)$, then Re $\psi(\lambda_0+(1-\lambda_0)u_2i,v_1)\leq\lambda_0$  if and only if $\mu_1+\lambda_0\geq 0$. This is true if Re $\mu=\mu_1\geq 0$ since $\lambda_0$ is nonegative. Thus by Lemma 1, if Re $\mu\geq 0$, then Re $\psi(p(z),zp'(z))>\lambda_0$ implies Re $p(z)>\lambda_0$. That is
\[
Re\;\frac{D^{n+1}\mathcal{J}_m^j(f)^\beta}{D^n\mathcal{J}_m^j(f)^\beta}>\lambda_0\;\;\text{if}\;\;Re\;\frac{D^{n+1}\mathcal{J}_{m-1}^j(f)^\beta}{D^n\mathcal{J}_{m-1}^j(f)^\beta}>\lambda_0.
\]
Since $\mathcal{J}_0^j(f)^\beta=f(z)^\alpha$ we have Re $\frac{D^{n+1}f(z)^\alpha}{D^nf(z)^\alpha}>\lambda_0\Rightarrow$ Re $\frac{D^{n+1}\mathcal{J}_1^j(f)^\beta}{D^n\mathcal{J}_1^j(f)^\beta}>\lambda_0\Rightarrow$ Re $\frac{D^{n+1}\mathcal{J}_2^j(f)^\beta}{D^n\mathcal{J}_2^j(f)^\beta}>\lambda_0\Rightarrow...$ and so on for all $m\in N$. By Lemma 2, we have: Re $\frac{D^{n+1}f(z)}{D^nf(z)}>\frac{\lambda_0}{\alpha}\Rightarrow$ Re $\frac{D^{n+1}\mathcal{J}_1^j(f)}{D^n\mathcal{J}_1^j(z)}>\frac{\lambda_0}{\beta}\Rightarrow$ Re $\frac{D^{n+1}\mathcal{J}_2^j(f)}{D^n\mathcal{J}_2^j(f)}>\frac{\lambda_0}{\beta}\Rightarrow...$ and so on for all $m\in N$. By setting $\lambda_0=\alpha\lambda$ we have Theorem 1.
\end{proof}

The next theorem will leads us to the best possible inclusion relations.

\begin{theorem}
Let $\alpha\geq 0$. Suppose for $\alpha>0$, the real number $\lambda$ is defined such that $0\leq\alpha\lambda<1$. If
\[
Re\;\frac{D^{n+1}\mathcal{J}_{m-1}^j(f)^\beta}{D^n\mathcal{J}_{m-1}^j(f)^\beta}>\alpha\lambda,\;\;then\;\;
\frac{D^{n+1}\mathcal{J}_m^j(f)^\beta}{D^n\mathcal{J}_m^j(f)^\beta}\prec q(z)
\]
where
\begin{equation}
q(z)=\frac{z^{1+\mu}(1-z)^{-2(1-\alpha\lambda)}}{\int_0^zt^\mu(1-t)^{-2(1-\alpha\lambda)}dt}-1\, \label{9}
\end{equation}
and $\mu=\gamma$ for $j=1$ and $\mu=\gamma+m-1$ for $j=2$.
\end{theorem}

\begin{proof}
As in the preceeding theorem, the case $\alpha=0$ holds trivially. Now for $\alpha>0$, let $0\leq\lambda_0=\alpha\lambda<1$ and suppose
\[
Re\;\frac{D^{n+1}\mathcal{J}_{m-1}^j(f)^\beta}{D^n\mathcal{J}_{m-1}^j(f)^\beta}>\lambda_0.
\]
Then from ~(\ref{8}), we have
\[
p(z)+\frac{zp'(z)}{\mu+p(z)}\prec\frac{1+(1-2\lambda_0)z}{1-z}
\]
Now by considering the differential equation
\[
q(z)+\frac{zq'(z)}{\mu+q(z)}=\frac{1+(1-2\lambda_0)z}{1-z}
\]
whose univalent solution is given by ~(\ref{9}) (see \cite{SS}), then by Lemma 2 we have the subordination
\[
p(z)=\frac{D^{n+1}\mathcal{J}_m^j(f)^\beta}{D^n\mathcal{J}_m^j(f)^\beta}\prec q(z)\prec\frac{1+(1-2\lambda_0)z}{1-z},
\]
where $q(z)$ is the best dominant, which proves the theorem.
\end{proof}

\section{The special case $m=1$}

\medskip
In this section we discuss the integral ~(\ref{2}), which coincides with the case $m=1$ of both integrals $\mathcal{J}_m^j(f)$. In particular we take $\lambda=0$. In this case, our first corollary, a simple one from Theorem 1, is the following:

\begin{corollary}
The classes $S_n$ is closed under $\mathcal{J}(f)$.
\end{corollary}

This result is more genral than the result of Miller et-al \cite{MM} (Theorem 2, pg. 162) in which case $\delta=\gamma$. A major breakthrough with our method is the fact that the integral ~(\ref{2}) passes through, preserving all the goemetry (starlikeness and convexity for example) of $f$ without having to drop any member of the sets on which the parameters $\alpha\geq 0$ and $\beta>0$ were defined, which was not the case in many earlier works. This will become more evident in the following more specific cases (cf. \cite{MM}).

{\rm (i)} If $\alpha+\delta=\beta+\gamma=1$, we have

\begin{corollary}
If $f\in S_n$, then
\[
\mathcal{J}(f)=\left\{z^{\beta-1}\int_0^z\left(\frac{f(t)}{t}\right)^\alpha dt\right\}^{\frac{1}{\beta}}=z+\cdots
\]
also belongs to $S_n$.
\end{corollary}

{\rm (ii)} If $\alpha+\delta=1$, $\beta=1$ and $\gamma=0$, we have

\begin{corollary}
If $f\in S_n$, then
\[
\mathcal{J}(f)=\int_0^z\left(\frac{f(t)}{t}\right)^\alpha dt=z+\cdots
\]
also belongs to $S_n$.
\end{corollary}

{\rm (iii)} If $\alpha+\delta=\beta+\gamma=\alpha+\eta+\gamma$, we have

\begin{corollary}
If $f\in S_n$, then
\[
\mathcal{J}(f)=\left\{\frac{\alpha+\gamma+\eta}{z^\gamma}\int_0^zt^{\gamma+\eta}f(t)^\alpha dt\right\}^{\frac{1}{\beta}}=z+\cdots
\]
also belongs to $S_n$.
\end{corollary}

From the above corollary, we can obtain various sequences of starlike and convex functions (and more generally, of $S_n$ functions): For example if $\gamma+\eta=1$, $\alpha=1$ and $\eta=k=0,1,2\cdots$; and if $\gamma=0$, $\alpha=1$ and $\eta=k=0,1,2\cdots$ we obtain, respectively, the following sequences of $S_n$ functions:
\[
\left\{2z^{k-1}\int_0^zf(t)dt\right\}^{\frac{1}{k+1}}=z+\cdots,\;\;\;k=0,\;1,\;2,\cdots
\]
and
\[
\left\{(k+1)\int_0^zt^{k-1}f(t)dt\right\}^{\frac{1}{k+1}}=z+\cdots,\;\;\;k=0,\;1,\;2,\cdots.
\]
For starlike functions, the above sequences are due to Miller et-al \cite{MM}.

Next we consider the best possible inclusion for the integral $\mathcal{J}(f)$ for two cases $\mu=0,1$. For these two cases, we have
\[
q(z)=\frac{1}{1-z},\;\;\;\mu=0,
\]
and
\[
q(z)=\frac{z^2}{(1-z)[(1-z)\ln(1-z)+z]}-1,\;\;\;\mu=1.
\]
But Re $q(z)\geq q(-r)$ for $|z|\leq r<1$. Thus we have
\[
Re\;q(z)\geq\frac{1}{1+r},\;\;\;\mu=0,
\]
and
\[
Re\;q(z)\geq\frac{1}{(1+r)[(1+r)\ln(1+r)-r]}-1,\;\;\;\mu=1.
\]

Letting $r\rightarrow 1^{-}$, we have Re $q(z)>\rho$ (say), $z\in E$. Since $\frac{D^{n+1}\mathcal{J}(f)^\beta}{D^n\mathcal{J}(f)^\beta}\prec q(z)$, we have Re $\frac{D^{n+1}\mathcal{J}(f)^\beta}{D^n\mathcal{J}(f)^\beta}>\rho$. By Lemma 3, this implies  Re $\frac{D^{n+1}\mathcal{J}(f)}{D^n\mathcal{J}(f)}>\frac{\rho}{\beta}$ so that the following best possible inclusions follow.

\begin{corollary}
Let $f\in S_n$. If $\delta$ is a real number and $\gamma=0$, then
\[
\mathcal{J}(f)=\left\{\beta\int_0^zt^{\delta-1}f(t)^\alpha dt\right\}^{\frac{1}{\beta}}=z+\cdots
\]
belongs to $S_n(\frac{1}{2\beta})$. 
\end{corollary}

\begin{corollary}
Let $f\in S_n$. If $\delta$ is a real number and $\gamma=1$, then
\[
\mathcal{J}(f)=\left\{\frac{\beta+1}{z}\int_0^zt^{\delta-1}f(t)^\alpha dt\right\}^{\frac{1}{\beta}}=z+\cdots
\]
belongs to $S_n(\frac{3-4\ln2}{2\beta(2\ln2-1)})$. 
\end{corollary}

On the final note, if we take $\beta=1$ in Corollaries 5 and 6 we then have the following special cases

\begin{corollary}
Let $f\in S_n$. If $\delta$ is a real number and $\gamma=0$, then
\[
\mathcal{J}(f)=\int_0^zt^{\delta-1}f(t)^\alpha dt=z+\cdots
\]
belongs to $S_n(\frac{1}{2})$. 
\end{corollary}

\begin{corollary}
Let $f\in S_n$. If $\delta$ is a real number and $\gamma=1$, then
\[
\mathcal{J}(f)=\frac{2}{z}\int_0^zt^{\delta-1}f(t)^\alpha dt=z+\cdots
\]
belongs to $S_n(\frac{3-4\ln2}{2(2\ln2-1)})$. 
\end{corollary}

Miller et-al \cite{MM} proved that if $f\in S^\ast$, then $\mathcal{J}(f)\in S^\ast(\frac{\sqrt{17}-3}{4})$ and also if $f\in K$, then $\mathcal{J}(f)\in K(\frac{\sqrt{17}-3}{4})$, which our last corollary has now raised to their best-possible status.
\medskip

{\it Acknowledgements.} This work was carried out at the Centre for Advanced Studies in Mathematics, CASM, Lahore University of Management Sciences, Lahore, Pakistan during the author's postdoctoral fellowship at the Centre. The author is indebted to all staff of CASM for their hospitality.

\vspace{10pt}

\hspace{-4mm}{\small{Received}}

\vspace{-12pt}
\ \hfill \
\begin{tabular}{c}
{\small\em  {\bf Current Address}}\\
{\small\em  Centre for Advanced Studies in Mathematics}\\
{\small\em  Lahore University of Management Sciences}\\
{\small\em  Lahore, Pakistan}\\
{\small\em E-mail: {\tt kobabalola@lums.edu.pk}}\\
{\small\em  {\bf Permanent Address}}\\
{\small\em  Department of Mathematics}\\
{\small\em  University of Ilorin}\\
{\small\em  Ilorin, Nigeria}\\
{\small\em E-mail: {\tt kobabalola@gmail.com}}\\
\end{tabular}

\end{document}